# On Knotted Subgroups

Marc Aaron F. Julian[1,*], Mark Lexter D. De Lara[1], Krizal John C. Espacio[1], Micko Jay S. Bajamundi[1], Clarisson Rizzie P. Canlubo[2]
[1]Institute of Mathematical Sciences, University of the Philippines Los Baños
[2]College of Computing and Information Technologies, National University−Manila
[*]Corresponding author: mfjulian@up.edu.ph

**Abstract**

In this article, we defined a knotted subgroup of a Lie group and considered a geometric notion of equivalence among them. We characterized these knotted subgroups in terms of one-parameter subgroups and provided examples in the case of $SU(2)$ and $SU(3)$. Infinitesimal elements that give rise to knotted subgroups of $SU(n)$ and $SO(n)$ are characterized as well. Canonical forms for their knotted subgroups are presented and their properties are described in terms of the spectrum of the corresponding infinitesimal elements. Finally, knotted subgroups of $SL(2,\mathbb{R})$ are completely classified using direct computation while knotted subgroups of $SL(3,\mathbb{R})$ are completely classified using Jordan canonical forms.

# 1 Introduction

The phenomenon of linking and winding of curves and surfaces are best handled using knots. In the most basic sense, a knot is a continuous mapping of $S^1$ to $\mathbb{R}^3$ which we denote by $S^1 \xrightarrow{\phi} \mathbb{R}^3$. We say that a knot $S^1 \xrightarrow{\phi} \mathbb{R}^3$ is *ambient isotopic* to another knot $S^1 \xrightarrow{\psi} \mathbb{R}^3$ if there is a continuous map $\mathbb{R}^3 \times [0,1] \xrightarrow{F} \mathbb{R}^3$ such that $F_0$ is the identity map, $F_t$ is a homeomorphism for all $t \in [0,1]$, and $F_1 \circ \phi = \psi$. Here, $\mathbb{R}^3 \xrightarrow{F_t} \mathbb{R}^3$ is the map given by $F_t(x) = F(x,t)$. Ambient isotopy is an equivalence relation on knots, hence, with only topology under consideration, it gives the appropriate notion of equivalence among knots.

In a more general setting, one can replace $\mathbb{R}^3$ with a $3-$manifold $M$. In fact, for compactness consideration, one often replaces $\mathbb{R}^3$ with its one-point compactification $S^3$. For higher dimensional linking and winding phenomenon, one uses higher dimensional spheres in place of $S^1$ with the target space being 2 dimensions higher [1].

Lie groups are among the most interesting class of manifolds. A *Lie group* is a manifold equipped with a group structure wherein the group product and inversion are smooth maps. For instance, the collection of $n \times n$ matrices with complex entries and nonzero determinant, denoted by $\mathrm{GL}_n(\mathbb{C})$, is a Lie group [2]. One can think of $\mathrm{GL}_n(\mathbb{C})$ as an open subset of $\mathbb{R}^{2n^2}$ since the determinant map $\mathbb{R}^{2n^2} \cong M_n(\mathbb{C}) \xrightarrow{\det} \mathbb{C}$ is smooth (being a polynomial) and $\mathrm{GL}_n(\mathbb{C}) = \det^{-1}(\mathbb{C}\backslash\{0\})$. The smooth structure on $\mathrm{GL}_n(\mathbb{C})$ is the one it inherits from $\mathbb{R}^{2n^2}$ as an open subset. The group operation is a smooth map since the entries of the product $AB$ are polynomials in terms of the entries of $A$ and $B$. The inversion map is also smooth since the inverse $A^{-1}$ of an invertible matrix $A$ is a polynomial in $A$ by the Cayley-Hamilton Theorem [3].

The Lie group $\mathrm{GL}_n(\mathbb{C})$ has a class of interesting closed subgroups called *matrix Lie groups.* Equivalently, a subgroup $G$ of $\mathrm{GL}_n(\mathbb{C})$ is a matrix Lie group if and only if any sequence of matrices $\{A_n\}$ in $G$ that converges to a matrix $A$, either has $A$ in $G$ or $A$ is non-invertible [4]. Let us consider examples of matrix Lie groups. Given a positive integer $n$, we denote by $U(n)$ the set of all matrices $X$ such that $X^* = X^{-1}$. We call this the *unitary group*. The subgroup $SU(n)$ of $U(n)$ which consists of unitary matrices with determinant 1 is called the *special unitary group*. We may also consider $\mathrm{GL}_n(\mathbb{R})$ as a matrix Lie group that sits inside $\mathrm{GL}_n(\mathbb{C})$ which also contains some interesting matrix Lie groups. For example, the group of all matrices $X$ satisfying $X^T = X^{-1}$, denoted by $O(n)$, is called the *orthogonal group*. Meanwhile, its subgroup which consists of orthogonal matrices with determinant 1, denoted by $SO(n)$ is called the *special orthogonal group*.

A Lie group $G$, being a manifold, naturally contains knots. Keeping in mind of its smooth structure means that we may focus our attention on smooth knots. The interesting question to answer would be: which among these (smooth) knots in $G$ inherit a group structure from $G$? We will call these *knotted subgroups* of $G$. More precisely, a subgroup $H$ of a Lie group $G$ is *knotted* if there is a continuous injective group homomorphism $S^1 \xrightarrow{\gamma} G$ with continuous inverse such that $\gamma(S^1) = H$. Oftentimes, we will refer to $S^1 \xrightarrow{\gamma} G$ as the knotted subgroup itself. Now, let us consider a notion of equivalence among knotted subgroups which must reflect the fact that knotted subgroups have a topology and a group structure that needs to be preserved. Let $S^1 \xrightarrow{\phi_1} G$ and $S^1 \xrightarrow{\phi_2} G$ be knotted subgroups of $G$ with images $H_1$, and $H_2$, respectively. A continuous map $F : G \times [0,1] \to G$ is an *ambient automorphism* from $\phi_1$ to $\phi_2$ if $F_0$ is the identity map on $G$; $F_t$ is a homeomorphism and a group automorphism on $G$ for all $t$; and $F_1 \circ \phi_1 = \phi_2$. Here, $F_t(A) := F(A,t)$. In this case, we say that $\phi_1$ (resp., $H_1$) is *ambient automorphic* to $\phi_2$ (resp., $H_2$). It is immediate to see that ambient automorphism is an equivalence relation among the knotted subgroups of $G$.

No work seems to be done with how subgroups of a Lie group are knotted. However, the concepts of knots and groups are deeply intertwined since several groups are associated with knots. We have the fundamental group of a knot complement [5], and the braid groups [6] to mention a few. Knots also found lots of connections to other areas of mathematics. For instance, parallel theories between knots and primes exists as outlined in [7]. Knot invariants come in many forms as well such as combinatorial invariants arising from projections of knots [8], quantum invariants from the theory of quantum groups and quantum algebras [6], and the Jones polynomials arising from operator algebras [9]. Knots also appear in fluid mechanics and dynamical systems as outlined in a talk by Etienne Ghys in the 2006 ICM. Knotted structures seem inevitable in mathematics and physics and this article is a reflection of this principle.

# 2 One-Parameter Subgroups as Knotted Subgroups

A *one-parameter subgroup* of a Lie group $G$ is a continuous group homomorphism $\mathbb{R} \xrightarrow{\phi} G$. Given a matrix Lie group $G$, there is a nice way to enumerate all one-parameter subgroups of $G$. If $\mathfrak{g}$ is the Lie algebra of $G$, then for any $X \in \mathfrak{g}$, the map

$$\begin{aligned} \phi_X : \ \mathbb{R} &\to G \\ t &\mapsto e^{tX} := \sum_{n=0}^{\infty} \frac{(tX)^n}{n!} \end{aligned} \tag{1}$$

is a one-parameter subgroup of $G$ [4]. In fact, any one-parameter subgroup of $G$ has this form. Furthermore, the one-parameter subgroup $\phi_X$ is trivial if and only if $X = 0$.

**Lemma 1.** *Let $G$ be a matrix Lie group with Lie algebra $\mathfrak{g}$. Elements $X, Y \in \mathfrak{g}$ are conjugate if and only the corresponding one-parameter subgroups are conjugate.*

*Proof:* Assume $X = P^{-1}YP$ for some appropriate $P$. Then,

$$\phi_X(t) = e^{tX} = e^{tP^{-1}YP} = \sum_{k=0}^{\infty} \frac{(tP^{-1}YP)^k}{k!} = P^{-1} \sum_{k=0}^{\infty} \frac{(tY)^k}{k!} P = P^{-1}e^{tY}P = P^{-1}\phi_Y(t)P.$$

Conversely, assume $\phi_X(t) = P^{-1}\phi_Y(t)P$ for some $P$. Then,

$$X = \left.\frac{d}{dt}\right|_{t=0} e^{tX} = \left.\frac{d}{dt}\right|_{t=0} \phi_X(t) = \left.\frac{d}{dt}\right|_{t=0} P^{-1}\phi_Y(t)P = \left.\frac{d}{dt}\right|_{t=0} e^{tP^{-1}YP} = P^{-1}YP.$$

This proves the lemma. ■

**Lemma 2.** *Any nontrivial proper subgroup of $\mathbb{R}$ is either infinite cyclic or dense in $\mathbb{R}$.*

**Theorem 3.** *There exists a one-to-one correspondence between knotted subgroups of $G$ and the nontrivial, non-injective one-parameter subgroups of $G$.*

*Proof.* Let $\phi : \mathbb{R} \longrightarrow G$ be a nontrivial and non-injective one-parameter subgroup of $G$. Then, Ker $\phi$ must be a nontrivial proper subgroup of $\mathbb{R}$. Since $\{I_n\}$ is closed in a matrix Lie group $G$, then Ker $\phi$ must be closed in $\mathbb{R}$. By Lemma 2, $\text{Ker}\phi = x\mathbb{Z}$ for some positive real number $x$. Without loss of generality, we may assume $\text{Ker}\phi = \mathbb{Z}$. By the First Isomorphism Theorem for Lie groups, $\phi(\mathbb{R}) \cong \mathbb{R}/\text{Ker}\phi = \mathbb{R}/\mathbb{Z} \cong S^1$. Thus, $\phi(\mathbb{R})$ is isomorphic to $S^1$ as a Lie group. This implies that there is an isomorphism $\gamma : S^1 \longrightarrow \phi(\mathbb{R}) \subseteq G$, making $\phi(\mathbb{R})$ a knotted subgroup of $G$.

Conversely, suppose $\gamma : S^1 \longrightarrow G$ is a knotted subgroup of $G$. Then $\gamma$ is a continuous group homomorphism with a continuous inverse. Consider the exponential map $\varphi : \mathbb{R} \longrightarrow S^1$ defined by $\varphi(x) = e^{ix}$, which is a continuous, nontrivial and non-injective group homomorphism from $\mathbb{R}$ to $S^1$. Consequently, the composition $\gamma \circ \varphi : \mathbb{R} \longrightarrow G$ is a nontrivial, non-injective one-parameter subgroup of $G$. ■

It is important to note that certain matrix Lie groups do not have knotted subgroups. For instance, consider the real Heisenberg group $H$ which consists of $3 \times 3$ real matrices of the form

$$\begin{pmatrix} 1 & a & c \\ & 1 & b \\ & & 1 \end{pmatrix}.$$

Its Lie algebra, denoted by $\mathfrak{h}$, consists of $3 \times 3$ real matrices of the form

$$\begin{pmatrix} & a & c \\ & & b \\ & & \end{pmatrix}.$$

In this case, the exponential map $\mathfrak{h} \xrightarrow{\exp} H$ is a bijection, thus, there is no $X \in \mathfrak{h}$ such that $\phi(t) = e^{tX}$ is non-injective.

# 3 Knotted Subgroups of $SU(n)$

Any element of $SU(2)$ can be written as

$$\begin{pmatrix} \alpha & -\bar{\beta} \\ \beta & \bar{\alpha} \end{pmatrix}$$

for complex numbers $\alpha$ and $\beta$ satisfying $|\alpha|^2 + |\beta|^2 = 1$. This, incidentally, gives a straightforward diffeomorphism between $SU(2)$ and the 3-dimensional sphere $S^3$. In particular, any element of $SU(2)$ can be written as $x_0 I_2 + x_1 \sigma_X + x_2 \sigma_Y + x_3 \sigma_Z$ where

$$\sigma_X = \begin{pmatrix} i & \\ & -i \end{pmatrix}, \quad \sigma_Y = \begin{pmatrix} & 1 \\ -1 & \end{pmatrix}, \quad \sigma_Z = \begin{pmatrix} & i \\ i & \end{pmatrix},$$

and $x_0^2 + x_1^2 + x_2^2 + x_3^2 = 1$ for $x_0, x_1, x_2, x_3 \in \mathbb{R}$. Consider the subset $Q \subseteq SU(2)$ consisting of those elements with $x_0 = 0$. Then, we have the following proposition whose proof is a straightforward computation.

**Proposition 4.** *Given any element $\Sigma \in Q$, the map $\mathbb{R} \xrightarrow{\phi_\Sigma} SU(2)$ given by $t \longmapsto (\cos t)\, I_2 + (\sin t)\, \Sigma$ is a one-parameter subgroup of $SU(2)$.*

Note that none of the maps $\phi_\Sigma$ is the trivial map and the periodicity of both the sine and cosine functions imply that these maps are non-injective. Therefore, these maps give rise to knotted subgroups. In particular, the knotted subgroups $\gamma_\Sigma$ they generate are the following

$$\begin{array}{rcl} S^1 & \xrightarrow{\gamma_\Sigma} & SU(2) \\ t (\text{mod } 2\pi) & \longmapsto & (\cos t)\, I_2 + (\sin t)\, \Sigma \end{array}$$

Next, let us consider the group $SU(3)$. The matrix Lie group $U(1)$ is isomorphic to $S^1$ as a Lie group. The group $U(1) \times U(1)$ is contained in $SU(3)$ using the following map

$$\begin{array}{rcl} U(1) \times U(1) & \longrightarrow & SU(3) \\ \left(e^{i\alpha}, e^{i\beta}\right) & \longmapsto & \begin{pmatrix} e^{i\alpha} & & \\ & e^{i\beta} & \\ & & e^{-i(\alpha+\beta)} \end{pmatrix} \end{array}$$

Since $U(1) \times U(1)$ is a torus of dimension 2, the image of the map above is a maximal torus of $SU(3)$. On the other hand, $U(1) \times U(1)$ can also be regarded as a surface in $\mathbb{R}^3$. For any coprime positive integers $p$ and $q$, the map $\gamma_{p,q}$ given by

$$\begin{array}{rcl} S^1 \cong \mathbb{R}/2\pi\mathbb{Z} & \longrightarrow & U(1) \times U(1) \subseteq \mathbb{R}^3 \\ t \ (\text{mod } 2\pi) & \longmapsto & (e^{ipt}, e^{iqt}) \end{array}$$

is a knot known as a *torus knot*. Figure 1 illustrates examples of torus knots. It is immediate to check that the map $\gamma_{p,q}$ is a group homomorphism. In particular, this implies that torus knots appear as knotted subgroups of $SU(3)$.

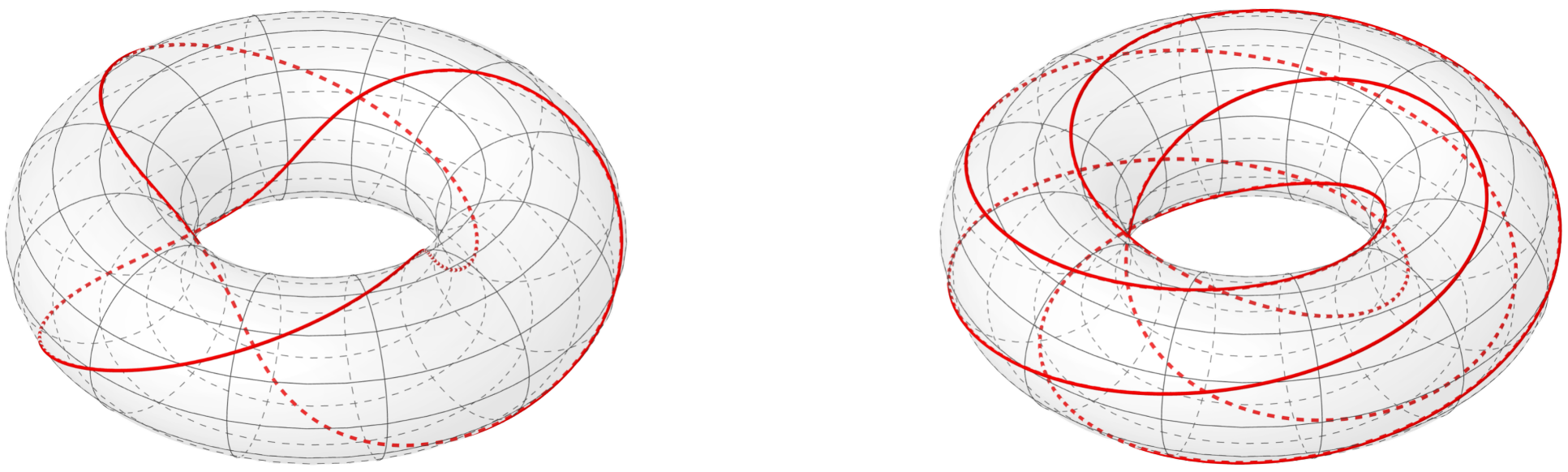

Figure 1: Torus knots: trefoil knot $\gamma_{3,2}$ (left) and $\gamma_{3,5}$ (right).

Finally, let us consider the groups $SU(n)$ for any $n \geqslant 4$. We denote by $E_{r,s}$ the $n \times n$ matrix with 1 at its $(r,s)-$entry and 0 elsewhere. The product of two matrices $E_{r,s}$ and $E_{v,w}$ is given by $E_{r,s}E_{v,w} = \delta_{s,v}E_{r,w}$ where $\delta_{s,v}$ is the Kronecker delta. Consider the following Pauli-type basis for $\mathfrak{su}(n)$

$$\begin{array}{rcllr} H_\ell & = & i(E_{\ell,\ell} - E_{\ell+1,\ell+1}), & & 1 \leq \ell \leq n-1 \\ X_{r,s} & = & E_{r,s} - E_{s,r}, & \text{where} & 1 \leq r < s \leq n \\ Y_{r,s} & = & i(E_{r,s} + E_{s,r}), & & 1 \leq r < s \leq n \end{array} . \tag{2}$$

The one-parameter subgroup generated by each basis element in 2 can be computed using equation 1 and are given by

$$\phi_{H_\ell}(t) = \left(\begin{array}{c:cc:c} I_{\ell-1} & & & \\ \hdashline & e^{it} & & \\ & & e^{-it} & \\ \hdashline & & & I_{n-\ell-1} \end{array}\right), \quad \phi_{X_{r,s}}(t) = \left(\begin{array}{c:c:c:c:c} I_{r-1} & & & & \\ \hdashline & \cos t & & \sin t & \\ \hdashline & & I_{s-r-1} & & \\ \hdashline & -\sin t & & \cos t & \\ \hdashline & & & & I_{n-s} \end{array}\right),$$

and

$$\phi_{Y_{r,s}}(t) = \left(\begin{array}{c:c:c:c:c} I_{r-1} & & & & \\ \hdashline & \cos t & & i\sin t & \\ \hdashline & & I_{s-r-1} & & \\ \hdashline & -i\sin t & & \cos t & \\ \hdashline & & & & I_{n-s} \end{array}\right).$$

Each of these are non-trivial and non-injective, hence, they give rise to knotted subgroups of $SU(n)$. Let us denote by $\gamma_{H_\ell}$, $\gamma_{X_{r,s}}$ and $\gamma_{Y_{r,s}}$ the knotted subgroups determined by $\phi_{H_\ell}$, $\phi_{X_{r,s}}$ and $\phi_{Y_{r,s}}$, respectively. The matrix representation of $\gamma_{H_\ell}(t)$ is the same as that of $\phi_{H_\ell}(t)$ but with $t$ taking values from $[0, 2\pi]$ instead of $\mathbb{R}$. Note that the matrix $\gamma_{H_\ell}(t)$ is permutation equivalent to $\gamma_{H_1}(t)$. Thus, there is a permutation matrix $P \in U(n)$ such that $\gamma_{H_1}(t) = P\gamma_{H_\ell}(t)P^{-1}$. Since $U(n)$ is path-connected, there is a path $[0,1] \xrightarrow{\xi} U(n)$ with $\xi(0) = I_n$ and $\xi(1) = P$. Define the map $SU(n) \times [0,1] \xrightarrow{F} SU(n)$ by $F(X,t) =$ $\xi(t)X\xi(t)^{-1}$ which we will also denote as $SU(n) \xrightarrow{F_t} SU(n)$ for $t \in [0,1]$. Note that this map is well-defined since for any $X \in SU(n)$ and $t \in [0,1]$, we have $F(x,t) \in SU(n)$.

It is immediate to see that $F_0$ is the identity map on $SU(n)$. For a fixed $t \in [0,1]$, $F_t$ acts as conjugation on $SU(n)$ by an element of $U(n)$, which is clearly a continuous automorphism with continuous inverse. Lastly,

$$F(\gamma_{H_\ell}(s), 1) = \xi(1)\gamma_{H_\ell}(s)\xi(1)^{-1} = P\gamma_{H_\ell}(s)P^{-1} = \gamma_{H_1}(s)$$

which implies that $F_1 \circ \gamma_{H_\ell} = \gamma_{H_1}$. Thus, $\gamma_{H_\ell}(t)$ is ambient automorphic to $\gamma_{H_\ell}(t)$.

Using the same argument, the knotted subgroups $\gamma_{X_{r,s}}$ are ambient automorphic to the knotted subgroup

$$\gamma_{X_{1,2}}(t) = \left(\begin{array}{cc|c} \cos t & \sin t & \\ -\sin t & \cos t & \\ \hline & & I_{n-2} \end{array}\right)$$

while the knotted subgroups $\gamma_{Y_{r,s}}$ are ambient automorphic to the knotted subgroup

$$\gamma_{Y_{1,2}}(t) = \left(\begin{array}{cc|c} \cos t & i\sin t & \\ -i\sin t & \cos t & \\ \hline & & I_{n-2} \end{array}\right).$$

Note that the knotted subgroup $\gamma_{X_{1,2}}$ is ambient automorphic to $\gamma_{H_1}$. To see this, consider the following special unitary matrix

$$P = \left(\begin{array}{cc|c} \frac{\sqrt{2}}{2} & \frac{\sqrt{2}}{2}i & \\ \frac{\sqrt{2}}{2}i & \frac{\sqrt{2}}{2} & \\ \hline & & I_{n-2} \end{array}\right).$$

Since $SU(n)$ is path-connected, there is a path $[0,1] \xrightarrow{\xi} SU(n)$ such that $\xi(0) = I_n$ and $\xi(1) = P$. Hence, the map $SU(n) \times [0,1] \xrightarrow{F} SU(n)$ defined as $F(X,t) = \xi(t)^{-1}X\xi(t)$ is an ambient automorphism from $\gamma_{X_{1,2}}$ to $\gamma_{H_1}$. Likewise, the knotted subgroup $\gamma_{Y_{1,2}}$ is ambient automorphic to $\gamma_{H_1}$.

**Theorem 5.** *For $n \geqslant 2$, any knotted subgroup of $SU(n)$ is ambient automorphic to a knotted subgroup of the form $\gamma(t) = \operatorname{diag}\left(e^{it\beta_1}, e^{it\beta_2}, \ldots, e^{it\beta_n}\right)$ where $\beta_\ell \in \mathbb{R}$.*

*Proof:* Consider an element $X \in \mathfrak{su}(n)$ such that the one-parameter subgroup generated by $X$, denoted by $\phi_X(t) = e^{tX}$, is non-trivial and non-injective. By the Spectral Theorem for $\mathfrak{su}(n)$, there is a unitary matrix $P$ such that $X = PDP^{-1}$ where $D$ is a diagonal matrix whose entries are the eigenvalues of $X$. Since $X$ is skew-hermitian, the eigenvalues of $X$ are purely imaginary numbers. Thus, the diagonal elements of $D$ are of the form $i\beta_\ell$ where $\beta_\ell \in \mathbb{R}$. Now,

$$\begin{aligned} P^{-1}\phi_X(t)P &= P^{-1}e^{tX}(t)P = P^{-1}\left(\sum_{k=0}^{\infty}\frac{(tX)^k}{k!}\right)P = \sum_{k=0}^{\infty}\frac{P^{-1}(tX)^k P}{k!} \\ &= \sum_{k=0}^{\infty}\frac{(tP^{-1}XP)^k}{k!} = \sum_{k=0}^{\infty}\frac{(tD)^k}{k!} = e^{tD} = \text{diag}(e^{it\beta_1}, e^{it\beta_2}, \ldots, e^{it\beta_n}). \end{aligned} \quad (3)$$

Since $U(n)$ is path-connected, there is a path $[0,1] \xrightarrow{\xi} U(n)$ such that $\xi(0) = I_n$ and $\xi(1) = P$. Hence, the map $SU(n) \times [0,1] \xrightarrow{F} SU(n)$ defined as $F(A,s) = \xi(s)^{-1}A\xi(s)$ is well-defined, i.e. for any $A \in SU(n)$ and $s \in [0,1]$ we have $F(A,s) \in SU(n)$. Therefore, $F$ is an ambient automorphism from $\phi_X$ to $\gamma(t) = e^{tD}$. This proves the theorem. ■

Theorem 5 implies that any knotted subgroup of SU($n$) is ambient automorphic to a knotted subgroup that is completely contained in a maximal torus[10]. Furthermore, it gives us a way to determine which elements of $\mathfrak{su}(n)$ give rise to knotted subgroups.

**Corollary 6.** *Let $X, Y \in \mathfrak{su}(n)$. If $X$ and $Y$ have the same list of eigenvalues then the knotted subgroup of $SU(n)$ generated by $X$ is ambient automorphic to the knotted subgroup generated by $Y$.*

**Corollary 7.** *Consider a non-zero $X \in \mathfrak{su}(n)$ whose eigenvalues are $i\beta_1, i\beta_2, \cdots, i\beta_n$, where $\beta_1, \beta_2, \cdots, \beta_n$ are real numbers. The one-parameter subgroup generated by $X$ is non-injective if and only if $\beta_1, \beta_2, \cdots, \beta_n$ are either all rational numbers or rational multiples of the same irrational number.*

# 4 Knotted Subgroups of $SO(n)$

An analogue of Theorem 5 holds true for the group $SO(n)$ but additional work is necessary since the proof of Theorem 5 does not readily translate to the case of $SO(n)$. In particular, infinitesimal elements of $SO(n)$ are not diagonalizable by orthogonal matrices. However, we can use a different canonical form for these infinitesimal elements.

Let $X \in \mathfrak{so}(n)$ which means that $X$ is a skew-symmetric matrix with real entries. By Schur Decomposition [3], there is a real orthogonal matrix $Q$, i.e. an element of $O(n)$, and a block diagonal matrix $D$ such that $X = Q^T DQ$. The block diagonal matrix $D$ consists of $1 \times 1$ and $2 \times 2$ blocks. The $1 \times 1$ blocks have entries equal to 0. The $2 \times 2$ blocks are of the form

$$J(\lambda) = \begin{pmatrix} & \lambda \\ -\lambda & \end{pmatrix}$$

where $\pm\lambda i$ is a conjugate pair of non-zero eigenvalues of $X$.

The group $O(n)$ has two connected components based on whether an element has determinant $+1$ or $-1$. Let us first consider the case when $Q$ has determinant $-1$. The orthogonal matrix $Q$ can then be factored as $Q = \Sigma\Omega$ where $\Omega \in SO(n)$ and

$$\Sigma = \left(\begin{array}{cc|c} & 1 & \\ 1 & & \\ \hline & & I_{n-2} \end{array}\right).$$

Then, $X = \Omega^T \Delta \Omega$ where $\Delta = \Sigma D \Sigma$. Let $\phi(t) = e^{tX}$ and $\psi(t) = e^{t\Delta}$ be the one-parameter subgroups determined by $X$ and $\Delta$, respectively. A similar computation as in (3), gives us

$$\begin{aligned} \Omega^T \psi(t) \Omega &= \Omega^T e^{t\Delta} \Omega = \Omega^T \left( \sum_{k=0}^{\infty} \frac{(t\Delta)^k}{k!} \right) \Omega \\ &= \sum_{k=0}^{\infty} \frac{\Omega^T (t\Delta)^k \Omega}{k!} = \sum_{k=0}^{\infty} \frac{\left(t\Omega^T \Delta \Omega\right)^k}{k!} = \sum_{k=0}^{\infty} \frac{(tX)^k}{k!} = e^{tX} = \phi(t). \end{aligned} \tag{4}$$

A straightforward computation shows that

$$\psi(t) = \left(\begin{array}{cccc|c} \lambda_1(t) & & & & \\ & \lambda_2(t) & & & \\ & & \ddots & & \\ & & & \lambda_\ell(t) & \\ \hline & & & & I_{n-2\ell} \end{array}\right)$$

where

$$\lambda_j(t) = \begin{pmatrix} \cos \lambda_j t & \sin \lambda_j t \\ -\sin \lambda_j t & \cos \lambda_j t \end{pmatrix}$$

and $\lambda_1, \dots, \lambda_\ell$ are non-zero real numbers. Thus, the one-parameter subgroup $\psi$ is non-trivial. If the sine and cosine functions here have synchronous periods then the one-parameter subgroup $\psi$ is also non-injective. In this case, $\psi$ determines a knotted subgroup of $SO(n)$. From equation (4), we see that $\phi$ is also non-trivial and non-injective.

Moreover, the knotted subgroup determined by $\phi$ is ambient automorphic to the knotted subgroup determined by $\psi$. By abuse of notation, let us denote by $\phi$ and $\psi$ these knotted subgroups, respectively. Since $SO(n)$ is path-connected, there is a path $[0,1] \xrightarrow{\xi} SO(n)$ such that $\xi(0) = I_n$ and $\xi(1) = \Omega$. Define the map $SO(n) \times [0,1] \xrightarrow{F} SO(n)$ by $F(A,s) = \xi(s)^T A \xi(s)$. Note that $F$ is continuous and well-defined, i.e. for any $A \in SO(n)$ and $s \in [0,1]$ we have $F(A,s) \in SO(n)$. By a direct computation, we see that $F$ is an ambient automorphism from $\psi$ to $\phi$.

The case in which $Q$ has determinant 1 is much easier since it is already an element of $SO(n)$ and the argument above works the same. Either way, we have the following.

**Theorem 8.** *Let $n \geqslant 2$. Any knotted subgroup of $SO(n)$ is ambient automorphic to a knotted subgroup of the form*

$$\gamma(t) = \left(\begin{array}{cccc|c} \lambda_1(t) & & & & \\ & \lambda_2(t) & & & \\ & & \ddots & & \\ & & & \lambda_\ell(t) & \\ \hline & & & & I_{n-2\ell} \end{array}\right)$$

*where* $\lambda_j(t) = \begin{pmatrix} \cos\lambda_j t & \sin\lambda_j t \\ -\sin\lambda_j t & \cos\lambda_j t \end{pmatrix}$ *and* $\lambda_1, \cdots, \lambda_\ell$ *are non-zero real numbers.*

An analogue of Corollary 6 holds true following Theorem 8 with $\mathfrak{so}(n)$ and $SO(n)$ in place of $\mathfrak{su}(n)$ and $SU(n)$, respectively. Similar to the case for $SU(n)$, Theorem 8 gives a characterization of those $X \in \mathfrak{so}(n)$ that give rise to knotted subgroups of $SO(n)$.

**Corollary 9.** *Consider a non-zero* $X \in \mathfrak{so}(n)$ *whose non-zero eigenvalues are* $\pm i\lambda_1$, $\pm i\lambda_2, \cdots, \pm i\lambda_\ell$. *The one-parameter subgroup generated by* $X$ *is non-injective if and only if* $\lambda_1, \lambda_2, \cdots, \lambda_\ell$ *are either all rational numbers or rational multiples of the same irrational number.*

# 5 Knotted Subgroups of $SL(2, \mathbb{R})$

The Lie algebra of $SL(2,\mathbb{R})$, denoted by $\mathfrak{sl}(2,\mathbb{R})$, consists of $2 \times 2$ traceless real matrices generated by

$$H = \begin{pmatrix} 1 & \\ & -1 \end{pmatrix}, \quad E = \begin{pmatrix} & 1 \\ 1 & \end{pmatrix}, \quad \text{and} \quad F = \begin{pmatrix} & -1 \\ 1 & \end{pmatrix}.$$

A direct computation of matrix exponentials using equation 1 show that the corresponding one-parameter subgroups are

$$\phi_H(t) = \begin{pmatrix} e^t & \\ & e^{-t} \end{pmatrix}, \quad \phi_E(t) = \begin{pmatrix} \cosh t & \sinh t \\ \sinh t & \cosh t \end{pmatrix}, \quad \text{and} \quad \phi_F(t) = \begin{pmatrix} \cos t & -\sin t \\ \sin t & \cos t \end{pmatrix},$$

respectively. Notice that $\phi_H$ and $\phi_E$ are injective since they both contain at least one injective entry. On the other hand, the one-parameter subgroup $\phi_F$ is non-injective. These are special instances of the following two propositions.

Let $a, b$ and $c$ be arbitrary real numbers. Consider the element $X = aE + bH + cF$ of $\mathfrak{sl}(2,\mathbb{R})$. Computing the corresponding one-parameter subgroup, we have the following

$$\phi(t) = e^{tX} = \frac{1}{\rho} \begin{pmatrix} \rho\cosh\rho t + b\sinh\rho t & (a-c)\sinh\rho t \\ (a+c)\sinh\rho t & \rho\cosh\rho t - b\sinh\rho t \end{pmatrix} \tag{5}$$

where $\rho = \sqrt{a^2+b^2-c^2}$. Note that if at least one of the entry functions of a matrix-valued function is injective, the matrix-valued function is injective. Using this principle, we have the following proposition.

**Proposition 10.** *Let $X = aE + bH + cF \in \mathfrak{sl}(2,\mathbb{R})$. If $\rho = \sqrt{a^2+b^2-c^2}$ is a non-zero real number then the one-parameter subgroup $\phi(t) = e^{tX}$ is injective.*

*Proof:* Assume $0 \neq \rho \in \mathbb{R}$. Let us consider cases depending on whether $a$ and $c$ are both 0 or not. If $a$ and $c$ are not both 0, then at least one of the $(1,2)$ or the $(2,1)$ entry functions are non-zero. In this case, $\phi$ has at least one injective entry function. If both $a$ and $c$ are 0, then $\rho = \pm b \neq 0$. Equation (5) then becomes

$$\phi(t) = e^{tX} = \begin{pmatrix} e^{\pm bt} & \\ & e^{\mp bt} \end{pmatrix},$$

which is clearly injective. In either case, we have $\phi$ injective. ■

**Proposition 11.** *Let $X = aE + bH + cF \in \mathfrak{sl}(2,\mathbb{R})$. The one-parameter subgroup $\phi(t) = e^{tX}$ is non-injective if $\rho = \sqrt{a^2+b^2-c^2}$ is a pure imaginary number.*

*Proof:* Suppose $\rho = i\lambda$ with $\lambda \neq 0$. Using the identities

$$\sinh(i\lambda t) = i\sin(\lambda t) \quad \text{and} \quad \cosh(i\lambda t) = \cos(\lambda t),$$

we can reduce the one-parameter subgroup in (5) into

$$\phi(t) = e^{tX} = \begin{pmatrix} \cos\lambda t + \frac{b}{\lambda}\sin\lambda t & \frac{a-c}{\lambda}\sin\lambda t \\ \frac{a+c}{\lambda}\sin\lambda t & \cos\lambda t - \frac{b}{\lambda}\sin\lambda t \end{pmatrix} \tag{6}$$

which is obviously periodic with period $\frac{2\pi}{\lambda}$. Thus, $\phi$ is non-injective. ■

**Proposition 12.** *Let $X = aE+bH+cF \in \mathfrak{sl}(2,\mathbb{R})$ be non-zero. If $\rho = \sqrt{a^2+b^2-c^2} = 0$, then $X$ is nilpotent of index 2 and the one-parameter subgroup $\phi(t) = e^{tX}$ is injective.*

*Proof:* A straightforward computation shows

$$X^2 = \begin{pmatrix} b & a-c \\ a+c & -b \end{pmatrix}^2 = \begin{pmatrix} b^2+a^2-c^2 & b(a-c)-b(a-c) \\ (a+c)b-(a+c)b & a^2-c^2+b^2 \end{pmatrix} = 0.$$

Thus, $X$ is nilpotent of index 2. Computing for the one-parameter subgroup gives us

$$e^{tX} = \begin{pmatrix} 1+bt & (a-c)\,t \\ (a+c)\,t & 1-bt \end{pmatrix}$$

which is clearly injective. ■

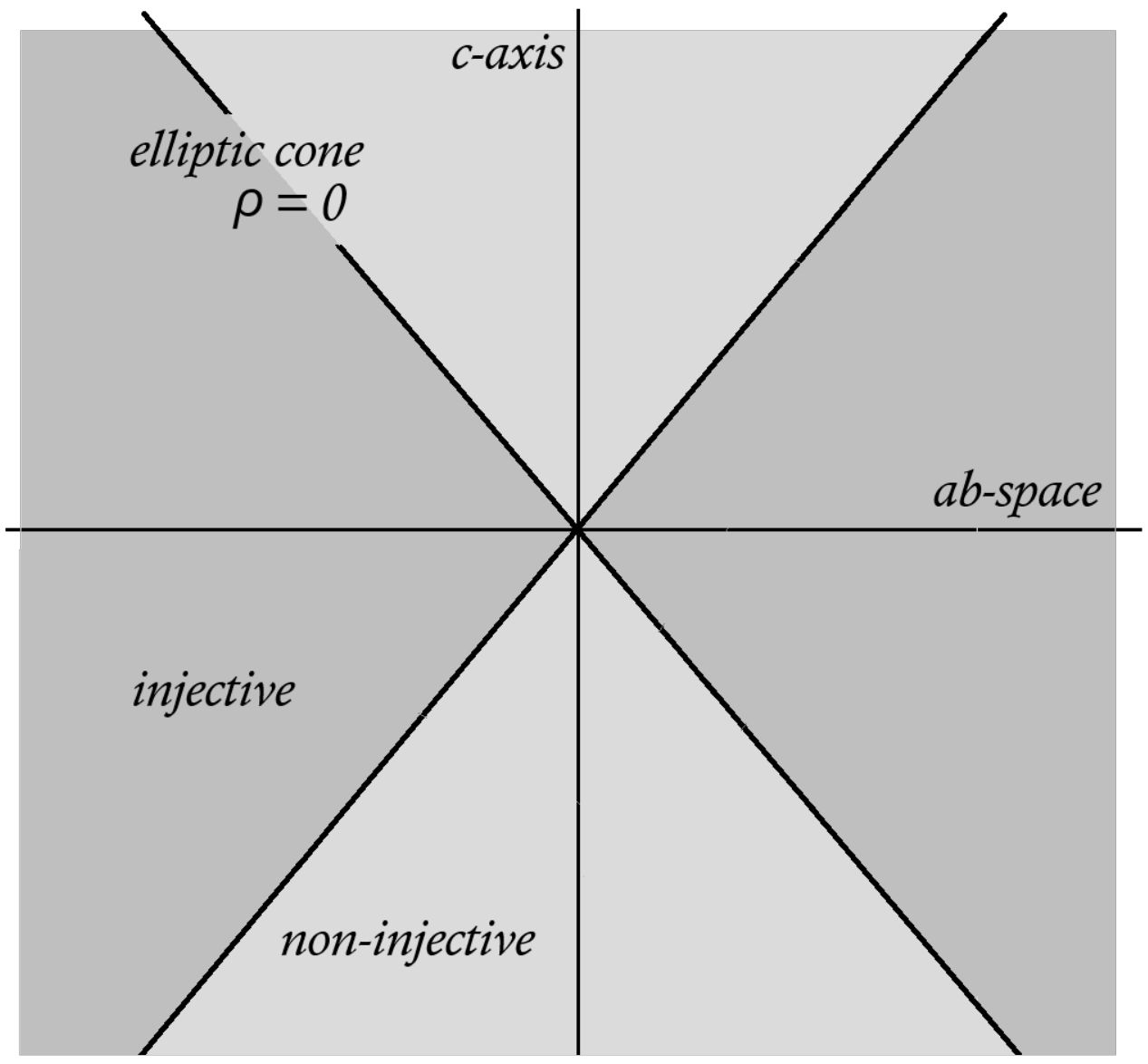

Figure 2: Parameter space of Theorem 13.

**Theorem 13.** *An element* $X = aE + bH + cF \in \mathfrak{sl}(2,\mathbb{R})$ *generates a knotted subgroup of* $SL(2,\mathbb{R})$ *if and only if* $a^2 + b^2 < c^2$.

There is a geometric way to visualize the current situation. Given an element $\left[\begin{smallmatrix} a & b \\ c & d \end{smallmatrix}\right]$ of $SL(2,\mathbb{R})$, where $ad - bc = 1$, the entries $a$ and $c$ cannot be simultaneously 0. Without loss of generality, assuming $a \neq 0$ means we can express $d = \frac{1+bc}{a}$. This implies that given three arbitrary real numbers where at least one is non-zero, a fourth real number is completely determined which will complete a $2 \times 2$ matrix element $A$ of $SL(2,\mathbb{R})$.

The first column of $A$ can be realized as points in the punctured plane $\mathbb{R}^2 \setminus (0,0)$. Then, for each such point of this puncture plane, we can choose a third arbitrary real number which completely determines an element of $SL(2,\mathbb{R})$. This implies that $SL(2,\mathbb{R})$ can be realized as a fibration over the puncture plane with fibers homeomorphic to $\mathbb{R}$ as shown in Figure 3. Finally, this space is homotopic to the infinite cylinder $S^1 \times \mathbb{R}$.

# 6 Knotted Subgroups of $SL(3,\mathbb{R})$

In Biswas et. al. [11], the one-parameter subgroups of $SL(3,\mathbb{R})$ are computed up to conjugacy. With additional work, we can compute the knotted subgroups of $SL(3,\mathbb{R})$ up to ambient automorphy. Consider the Lie algebra $\mathfrak{sl}(3,\mathbb{R})$ of $SL(3,\mathbb{R})$ whose elements are $3 \times 3$ real, traceless matrices. By the Jordan Decomposition Theorem, any element $X$ of $\mathfrak{sl}(3,\mathbb{R})$ with 3 real eigenvalues is conjugate to exactly one of the following forms

$$X_1 = \begin{pmatrix} \lambda_1 & & \\ & \lambda_2 & \\ & & -(\lambda_1 + \lambda_2) \end{pmatrix}, \quad X_2 = \begin{pmatrix} \lambda & 1 & \\ & \lambda & \\ & & -2\lambda \end{pmatrix}, \quad \text{and } X_3 = \begin{pmatrix} & 1 & \\ & & 1 \\ & & \end{pmatrix}.$$

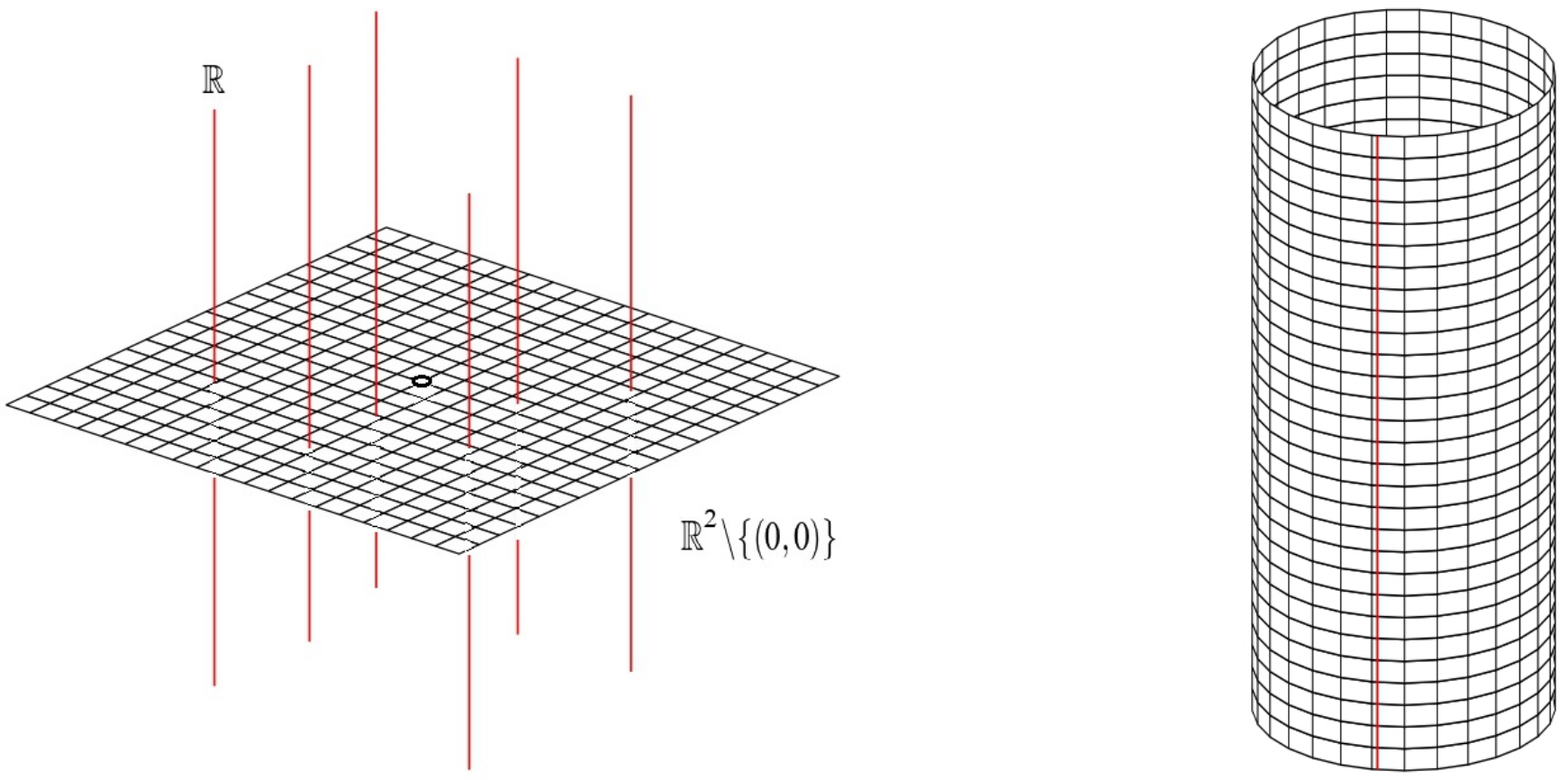

Figure 3: Fibration over $\mathbb{R}^2\backslash\{(0,0)\}$ with fibers $\mathbb{R}$ (left) and the infinite cylinder $S^1 \times \mathbb{R}$ (right).

If the eigenvalues of $X$ are distinct, it is conjugate to $X_1$. If $X$ has an eigenvalue with multiplicity 2, it is conjugate to $X_2$. If the eigenvalue of $X$ has multiplicity 3, it is conjugate to $X_3$. Now, any element $X$ of $\mathfrak{sl}(3,\mathbb{R})$ with a complex-conjugate pair of eigenvalues $a \pm bi$ is conjugate to an element of the following form

$$X_4 = \begin{pmatrix} a & b & \\ -b & a & \\ & & -2a \end{pmatrix}.$$

The corresponding one-parameter subgroups are

$$\phi_1(t) = e^{tX_1} = \begin{pmatrix} e^{\lambda_1 t} & & \\ & e^{\lambda_2 t} & \\ & & e^{-(\lambda_1+\lambda_2)t} \end{pmatrix}, \quad \phi_2(t) = e^{tX_2} = \begin{pmatrix} e^{\lambda t} & te^{\lambda t} & \\ & e^{\lambda t} & \\ & & e^{-2\lambda t} \end{pmatrix},$$

$$\phi_3(t) = e^{tX_3} = \begin{pmatrix} 1 & t & t^2/2 \\ & 1 & t \\ & & 1 \end{pmatrix}, \quad \text{and } \phi_4(t) = e^{tX_4} = \begin{pmatrix} e^{at}\cos bt & e^{at}\sin bt & \\ -e^{at}\sin bt & e^{at}\cos bt & \\ & & e^{-2at} \end{pmatrix}.$$

As we have seen in previous sections, conjugacy of elements of the Lie algebra lifts to conjugacy of the corresponding one-parameter subgroups. In particular, $\phi_1$ and all the one-parameter subgroups conjugate to $\phi_1$ are injective. This is because the $(1,1)$-entry function $e^{\lambda_1 t}$ is injective. Likewise, $\phi_2$ and $\phi_3$ and all the one-parameter subgroups conjugate to $\phi_2$ and $\phi_3$ are injective.

Notice that $\phi_4$ is injective if and only if $a \neq 0$. Thus, only the elements of $\mathfrak{sl}(3,\mathbb{R})$ with a pure imaginary eigenvalue give rise to a non-injective one-parameter subgroup.

**Theorem 14.** *Let $X \in \mathfrak{sl}(3,\mathbb{R})$. If $X$ has an eigenvalue $i\beta$ for some $0 \neq \beta \in \mathbb{R}$ then the one-parameter subgroup generated by $X$ is ambient automorphic to a one-parameter subgroup of the form*

$$\phi(t) = \begin{pmatrix} \cos\beta t & \sin\beta t & \\ -\sin\beta t & \cos\beta t & \\ & & 1 \end{pmatrix}.$$

*Proof.* That $X$ has a non-zero pure imaginary eigenvalue if and only if the one-parameter subgroup generated by $X$ is conjugate to $\phi$, follows immediately from the argument above and Lemma 1. All that is left to show is that the conjugacy of one-parameter subgroups lifts to ambient automorphy. Let $\psi(t) = e^{tX}$ be a one-parameter subgroup conjugate to $\phi$. Lemma 1 implies that there is a $P \in GL(3,\mathbb{R})$ satisfying $X = P^{-1}\Omega P$ where

$$\Omega = \begin{pmatrix} & \beta & \\ -\beta & & \\ & & \end{pmatrix}.$$

Note that $GL(3,\mathbb{R})$ has two connected components $GL_{\pm}(3,\mathbb{R})$ according to whether an element has a positive or a negative determinant.

Let us first consider the case when $P$ has negative determinant. In this case, we can factor $P$ as $P = \Sigma Q$ where $Q \in GL_{+}(3,\mathbb{R})$ and

$$\Sigma = \begin{pmatrix} & 1 & \\ 1 & & \\ & & 1 \end{pmatrix}.$$

Then, $X = P^{-1}\Omega P = Q^{-1}\Sigma^{-1}\Omega\Sigma Q$. Let $\Delta = \Sigma^{-1}\Omega\Sigma$. Then,

$$\Delta = \begin{pmatrix} & -\beta & \\ \beta & & \\ & & \end{pmatrix}.$$

Thus, $X = Q^{-1}\Delta Q$. Now, $GL_{+}(3,\mathbb{R})$ is path-connected. This implies that there is a path $[0,1] \xrightarrow{\xi} GL_{+}(3,\mathbb{R})$ such that $\xi(0) = I_3$ and $\xi(1) = Q$. Consider the map $SL(3,\mathbb{R}) \times [0,1] \xrightarrow{F} SL(3,\mathbb{R})$ given as $F(A,s) = \xi(s)^{-1}A\xi(s)$. Then, $F$ is well-defined, i.e. for any $A \in SL(3,\mathbb{R})$ and $s \in [0,1]$ we have $F(A,s) \in SL(3,\mathbb{R})$. Now, $F$ is an ambient automorphism from the one-parameter subgroup generated by $\Delta$ to the one-parameter subgroup $\psi$. Note that upon considering $-\beta$ instead of $\beta$, the one-parameter subgroup generated by $\Delta$ is $\phi$.

For the case when $P$ has positive determinant, a path $\xi$ from $I_3$ to $P$ readily exists and the construction of an ambient automorphism $F$ from $\phi$ to $\psi$ is much more straightforward. This proves the theorem. ∎

**Corollary 15.** *An element $X \in \mathfrak{sl}(3,\mathbb{R})$ generates a knotted subgroup if and only if $X$ has a pure imaginary eigenvalue and one zero eigenvalue.*

# 7 Epilogue

The Lie group $SU(n)$ has dimension $n^2 - 1$. In particular, $SU(2)$ has dimension 3 which guarantees the abundance of usual knots and their complexity. On the other hand, the dimension of $SU(n)$ for $n \geqslant 3$ is at least 8 which means that nothing interesting happens as far as knots are concerned [12]. However, a different story arise when dealing with knotted subgroups. The group structure heavily restricts the fibers that may appear in an ambient automorphism. Another subtlety that is present in the business of knotted subgroups is the trivial knot. In a simply-connected space without any regard to a group structure, the unknot is null-homotopic, i.e. it is homotopic to a point. For knotted subgroups, it does not make sense to speak of homotopy that preserves the group structure. This is because the unknot is isomorphic to the circle group $S^1$ while a point is the trivial group. Figure 4 shows the difference between a null-homotopic and a non-null-homotopic trivial knots.

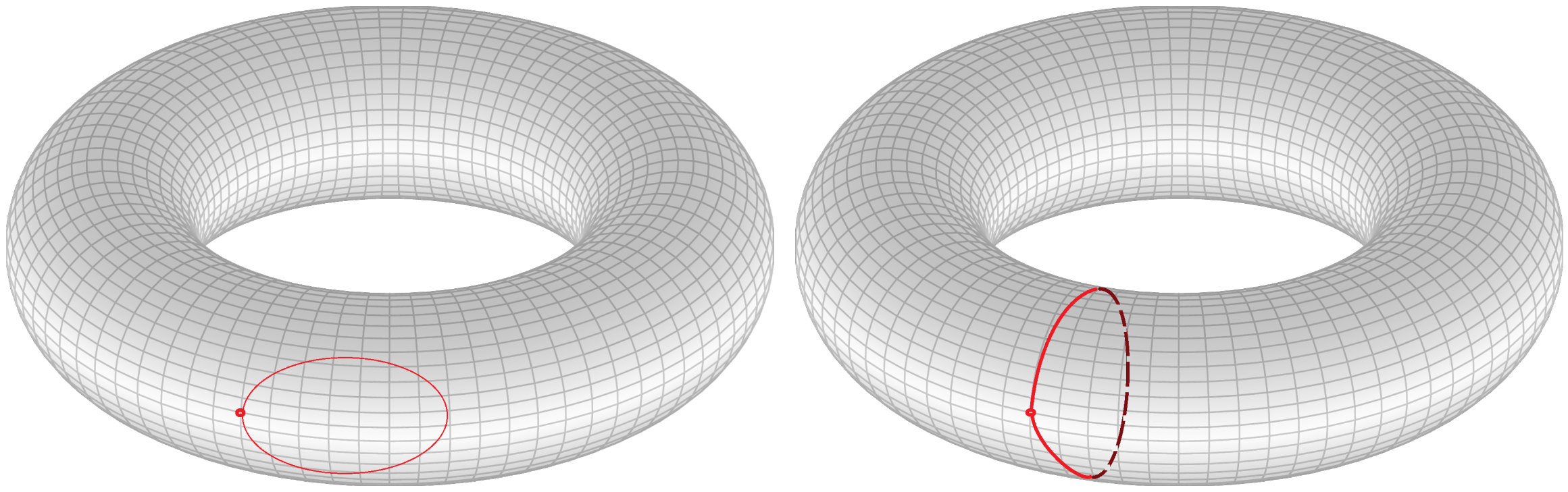

Figure 4: Null-homotopic (left) and non-null-homotopic (right) trivial knots.

Possible extension of this study is to investigate the knotted subgroups of other matrix Lie groups such as $SL(n, \mathbb{R})$ or $\mathrm{Sp}(n)$. The computations done in Sections 5 and 6 gives us a hint of which infinitesimal elements of $SL(n, \mathbb{R})$ generate a knotted subgroup. In particular, we have the following.

**Claim 16.** *An element $X \in \mathfrak{sl}(n, \mathbb{R})$ generates a knotted subgroup of $SL(n, \mathbb{R})$ if $X$ is diagonalizable and the eigenvalues of $X$ are pure imaginary numbers, not all of which are* 0*. In particular, the rank of $X$ is even.*

Preliminary results from Sections 5 and 6 suggest the following.

**Claim 17.** *The converse of Proposition 16 holds true, i.e. only diagonalizable elements of $\mathfrak{sl}(n, \mathbb{R})$ with pure imaginary eigenvalues, not all of which are* 0*, generate knotted subgroups of $SL(n, \mathbb{R})$.*

For $n = 2, 3$, it is also a curious observation that the infinitesimal elements $X$ of $SL(n, \mathbb{R})$ generating knotted subgroups are conjugate to elements of the proper Lie subalgebra $\mathfrak{so}(n)$ of $\mathfrak{sl}(n, \mathbb{R})$. This may be true for all $n \geqslant 2$, i.e.

**Claim 18.** *Let $n \geqslant 2$. Then, $X \in \mathfrak{sl}(n, \mathbb{R})$ generates a knotted subgroup of $SL(n, \mathbb{R})$ if and only if there is a $P \in GL(n, \mathbb{R})$ such that $P^{-1}XP \in \mathfrak{so}(n)$.*

Compactness of the subgroup $SO(n)$ of $SL(n, \mathbb{R})$ seems to be the reason for this but a thorough look is necessary.

Lastly, the most interesting aspect of knotted subgroups would be in the set up of 3-dimensional Lie groups since this is where the topology of these subgroups play a more crucial role at par with the group structure. It might be interesting to look at the ramifications of this study in relation to the Bianchi classification of 3-dimensional real Lie algebras.

## 8 Acknowledgment

This study is funded by the Collaborative Research Program of the Institute of Mathematical Sciences, University of the Philippines Los Baños.